\documentclass[11pt,a4paper]{article}
   \usepackage{amssymb}  \usepackage{amsmath}
   \usepackage{amsthm}   \usepackage{rotate}
   \usepackage{graphicx} \usepackage[T1]{fontenc}

\theoremstyle{definition}
\newtheorem{definition}{Definition}[section]
\theoremstyle{plain}
\newtheorem{theorem}[definition]{Theorem}
\newtheorem{lemma}[definition]{Lemma}

\newtheorem{corollary}[definition]{Corollary}

\begin{document}

\title{\bf On intuitionistic fuzzy sub-hyperquasigroups of hyperquasigroups}

\author{\bf Wies{\l}aw A. Dudek$^{\rm a}$, Bijan Davvaz$^{\rm b},$
             Young Bae Jun$^{\rm c,*}$ \\[2mm]
 {\small \it $^{\rm a}$ Institute of Mathematics, Technical University,} \\
 {\small \it Wybrze\.ze Wyspia\'nskiego 27, 50-370 Wroc{\l}aw, Poland}\\
 {\small \it $^{\rm b}$ Department of Mathematics, Yazd University, Yazd, Iran}\\
 {\small \it $^{\rm c}$ Department of Mathematics Educations,
               Gyeongsang National University,} \\
 {\small \it   Chinju 660-701, Korea}}

\date{}
 \maketitle

\begin{abstract}
The notion of intuitionistic fuzzy sets was introduced by
Atanassov as a generalization of the notion of fuzzy sets. In this
paper, we consider the intuitionistic fuzzification of the concept
of sub-hyperquasigroups in a hyperquasigroup and investigate some
properties of such sub-hyperquasigroups. In particular, we
investigate some natural equivalence relations on the set of all
intuitionistic fuzzy sub-hyperquasigroups of a hyperquasigroup.
\\[3mm]
 {\it 2000  Mathematics Subject Classification:} 20N20, 20N25. \\
 {\it Keywords:} hyperquasigroup, fuzzy sub-hyperquasigroup,
intuitionistic fuzzy sub-hyperquasigroup, quasigroup.
\end{abstract}

 \footnote{{\it E-mail address:}
        dudek\symbol{64}im.pwr.wroc.pl (W. A. Dudek),
        davvaz\symbol{64}yazduni.ac.ir (B. Davvaz),
        ybjun\symbol{64}nongae.gsnu.ac.kr (Y. B. Jun)
        (\tiny{hyper/His/INS394.tex})}

\section{Introduction and preliminaries}

The theory of hyperstructures which is a generalization of the
concept of algebraic structures first was introduced by Marty
\cite{14} and then many researchers have been worked on this new
field of modern algebra and developed it. A short review of the
theory of hyperstructures appear in \cite{4} and \cite{16}. A
recent book \cite{3} contains a wealth of applications. There are
applications to the following subjects: geometry, hypergraphs,
binary relations, lattices, fuzzy sets and rough sets, automata,
cryptography, combinatorics, codes, artificial intelligence, and
probabilities. The theory of fuzzy sets proposed by Zadeh
\cite{17} has achieved a great success in various fields. Out of
several higher order fuzzy sets, intuitionistic fuzzy sets
introduced by Atanassov \cite{1, 2, 2a} have been found to be
highly useful to deal with vagueness. Gau and Buehrer
\cite{IEEE23-610} presented the concept of vague sets. But,
Burillo and Bustince \cite{FSS79-403} showed that the notion of
vague sets coincides with that of intuitionistic fuzzy sets.
Szmidt and Kacprzyk \cite{FSS118-467} proposed a
non-probabilistic-type entropy measure for intuitionistic fuzzy
sets. De et al. \cite{FSS117-209} studied the Sanchez's approach
for medical diagnosis and extended this concept with the notion of
intuitionistic fuzzy set theory. Dengfeng and Chuntian
\cite{PRL23-221} introduced the concept of the degree of
similarity between intuitionistic fuzzy sets, presented several
new similarity measures for measuring the degree of similarity
between intuitionistic fuzzy sets, which may be finite or
continuous, and gave corresponding proofs of these similariry
measures and discussed applications of the similarity measures
between intuitionistic fuzzy sets to pattern recofnition problems.
The notion of join space has been introduced by Prenowitz and used
by him and afterwards together Jantosciak to build again several
branches of geometry. A join space is a hypergroup with additional
conditions. A generalization of join spaces for the point of view
of independence, dimension etc., is that of cambiste hypergroups
studied by Freni. Noticing that a hypergroup is a hyperquasigroup
with the associative hyperoperation, the results of this paper
will make a contribution to discuss a generalization of join
spaces, to deal with several notions in geometries since there are
deep relations between geometries and hypergroups (or, to say
multigroups), and to develop the intuitionistic fuzzy theory in
several algebraic structures.

A {\it hypergroupoid} $(G,\circ )$ is a non-empty set $G$ with a
{\it hyperoperation} $\circ$ defined on $G$, i.e., a mapping of
$G\times G$ into the family of non-empty subsets of $G$. If $(x,y)
\in G\times G,$ its image under $\circ$ is denoted by $x\circ y$.
If $A,B \subseteq G$ then $A\circ B$ is given by $A\circ B=\bigcup
\{x\circ y \ | \ x\in A, \ y\in B\}$. $x\circ A$ is used for $\{
x\}\circ A$ and $A\circ x$ for $A\circ \{ x\}$.

\begin{definition}\label{def11} A hypergroupoid $(G,\circ )$ is
called a {\it hypergroup} if for all $x,y,z\in G$ the following
two conditions hold:
\begin{itemize}
 \item[\rm (i)]  $x\circ (y\circ z)=(x\circ y)\circ z$,
 \item[\rm (ii)] $x\circ G = G\circ x=G$.
\end{itemize}
\end{definition}
The second condition, called the {\it reproduciblity condition},
means that for any $x,y\in G$ there exist $u,v\in G$ such that
$y\in x\circ u$ and $y\in v\circ x$.

A hypergroupoid satisfying this condition is called a {\it
hyperquasigroup}. Thus a hypergroup is a hyperquasigroup with the
associative hyperoperation.

A non-empty subset $K$ of a hyperquasigroup $(G,\circ )$ is called
a {\it sub-hyperquasigroup} if $(K,\circ)$ is a hyperquasigroup.

The concept of fuzzy sets was introduced by Zadeh \cite{17} in
1965. A mapping $\mu :X \to  [0,1]$, where $X$ is an arbitrary
non-empty set, is called a {\it fuzzy set} in $X$. The {\it
complement} of $\mu$, denoted by $\mu^c$, is the fuzzy set in $X$
given by $\mu^c(x)=1-\mu (x)$ for all $x \in X$.

For any fuzzy set $\mu$ in $X$ and any $t\in [0,1]$ we define two
sets
 \[U(\mu;t)=\{x\in X\ |\ \mu(x)\geq t\} \ \ \ \ {\rm and } \ \ \ \
   L(\mu;t)=\{x\in X\ |\ \mu(x)\leq t\}, \]
which are called an {\it upper } and {\it lower $t$-level cut} of
$\mu$ and can be used to the characterization of $\mu$.

 In 1971, Rosenfeld \cite{15} applied the concept of fuzzy sets to the
theory of groups and studied fuzzy subgroups of a group. Davvaz
applied in \cite{6} fuzzy sets to the theory of algebraic
hyperstructures and studied their fundamental properties. Further
investigations are contained in \cite{5}, \cite{7} and \cite{8}.

\begin{definition}(cf. \cite{6}) Let $(G,\circ )$ be a
hypergroup (resp. hyperquasigroup) and let $\mu$ be a fuzzy set in
$G$. Then $\mu$ is said to be a {\it fuzzy sub-hypergroup} (resp.
{\it fuzzy sub-hyperquasigroup}) of $G$ if the following axioms
hold:
\begin{itemize}
 \item[(1)] $\min\{\mu(x),\mu(y)\} \leq\inf\{\mu(z)\ |\ z\in x\circ
          y\}$ \ for all $x,y\in G$,
 \item[(2)] for all $x,a\in G$ there exists $y\in G$ such that $x\in a\circ y$ and
   \[ \min\{\mu(a),\mu(x)\}\leq\mu(y), \]
 \item[(3)] for all $x,a\in G$ there exists $z\in G$ such that $x\in z\circ a$ and
   \[\min\{\mu(a),\mu(x)\}\leq\mu(z).\]
\end{itemize}
\end{definition}

As an important generalization of the notion of fuzzy sets in $X$,
Atanassov \cite{1} introduced the concept of {\it intuitionistic
fuzzy sets} defined on a non-empty set $X$ as objects having the
form
 \[ A=\{(x,\mu_A(x), \lambda_A(x)) \ | \ x\in X\}, \]
where the functions $\mu_A:X\to  [0,1]$ and $\lambda_A:X\to
[0,1]$ denote the {\it degree of membership} (namely $\mu_A(x)$)
and the {\it degree of nonmembership} (namely $\lambda_A(x)$) of
each element $x\in X$ to the set $A$ respectively, and
$0\leq\mu_A(x)+\lambda_A(x) \leq 1$ for all $x \in X$.

Such defined objects are studied by many authors (see for example
two journals: 1. {\it Fuzzy Sets and Systems} and 2. {\it Notes on
Intuitionistic Fuzzy Sets}) and have many interesting applications
not only in mathematics (see Chapter 5 in the book \cite{2a}). In
particular, Kim, Dudek and Jun in \cite{11} introduced the notion
of an intuitionistic fuzzy subquasigroup of a quasigroup. Also in
\cite{12}, Kim and Jun introduced the concept of intuitionistic
fuzzy ideals of semigroups.

For every two intuitionistic fuzzy sets $A$ and $B$ in $X$ we
define (cf. \cite{2}):
\begin{itemize}
 \item[(1)] $A\subseteq B$ \ iff \ $\mu_A(x)\leq\mu_B(x)$ and
         $\lambda_A(x)\geq\lambda_B(x)$ \ for all $x\in X$,
 \item[(2)] $A^c=\{(x, \lambda_A(x),\mu_A (x)) \ | \ x\in X\}$,
 \item[(3)] $A\cap B=\{(x,\min\{\mu_A(x),\mu_B(x)\},
  \max\{\lambda_A(x),\lambda_B(x)\}) \ | \ x\in X\}$,
 \item[(4)] $A\cup B=\{(x,\max\{\mu_A(x),\mu_B(x)\},
            \min\{\lambda_A(x),\lambda_B(x)\}) \ | \ x\in X\}$,
 \item[(5)] $\Box A=\{(x,\mu_A(x),\mu^c_A(x)) \ | \ x\in X\}$,
 \item[(6)] $\Diamond A=\{(x,\lambda^c_A(x),\lambda_A(x)) \ | \ x\in X\}$.
\end{itemize}

\section{Intuitionistic fuzzy sub-hyperquasigroups}

For the sake of simplicity, we shall use the symbol $A=(\mu_A,
\lambda_A)$ for the intuitionistic fuzzy set
$A=\{(x,\mu_A(x),\lambda_A(x) \ | \ x\in X\}$.

In what follows, let $G$ denote a hyperquasigroup, and we start by
defining the notion of intuitionistic fuzzy sub-hyperquasigroups.

Based on \cite{11}, we can extend the concept of the
intuitionistic fuzzy subquasigroup to the concept of
intuitionistic fuzzy sub-hyperquasigroups in the following way:

\begin{definition}\label{def21} An intuitionistic fuzzy set $A=(\mu_A,
\lambda_A)$ in $G$ is called an {\it intuitionistic fuzzy
sub-hyperquasigroup} of $G$ ($IFSH$ of $G$ for short) if
\begin{itemize}
 \item[(1)] $\min\{\mu_A(x),\mu_A(y)\}\leq\inf\{\mu_A(z)\ |\
                z\in x\circ y\}$ \ for all $x,y\in G$,
 \item[(2)] for all $x,a\in G$ there exist $y,z\in G$ such that
$x\in (a\circ y)\cap (z\circ a)$ and
 \[ \min\{\mu_A(a),\mu_A(x)\}\leq\min\{\mu_A(y),\mu_A(z)\}, \]
 \item[(3)] $\sup\{\lambda_A(z)\ |\ z\in x\circ y\}\leq\max\{\lambda_A(x),\lambda_A(y)\}$
               \ for all $x,y\in G$,
 \item[(4)] for all $x,a\in G$ there exist $y,z\in G$ such that
  $x\in (a\circ y)\cap (z\circ a)$ and
 \[
 \max\{\lambda_A(y),\lambda_A(z)\}\leq\max\{\lambda_A(a),\lambda_A(x)\}.
 \]
\end{itemize}
\end{definition}
\begin{lemma}\label{lem22} If $A=(\mu_A, \lambda_A)$ is an
$IFSH$ of $G$, then so is \ $\Box A=(\mu_A, \mu^c_A)$.
\end{lemma}

\begin{proof} It is sufficient to show that $\mu_A^c$ satisfies the third and
fourth conditions of Definition \ref{def21}. For $x,y\in G$ we
have
$$
\min\{\mu_A(x),\mu_A(y)\}\leq\inf\{\mu_A(z)\ |\ z\in x\circ y\}
$$
and so
$$
\min\{1-\mu_A^c (x), 1-\mu_A^c (y)\} \leq\inf\{1-\mu_A^c (z)\ |\ z\in x\circ y\}.
$$
Hence
$$
\min\{1-\mu_A^c(x),1-\mu_A^c(y)\}\leq 1-\sup\{\mu_A^c(z)\ |\ z\in
x\circ y\}
$$
which implies
$$
\sup\{\mu_A^c (z)\ |\ z\in x\circ y\}\leq 1-\min\{1-\mu_A^c(x),
1-\mu_A^c(y)\}.
$$
Therefore
 $$ \sup\{\mu_A^c(z)\ | \ z\in x\circ y\}\leq\max\{\mu_A^c(x),\mu_A^c(y)\}. $$
Hence the third condition of Definition \ref{def21} is verified.

Now, let $a,x\in G.$ Then there exist $y,z\in G$ such that $x\in
a\circ y, \ x\in z\circ a$ and
\[
\min\{\mu_A (a),\mu_A (x)\} \leq \min\{\mu_A (y),\mu_A (z)\}.
\]
So
\[
\min\{1-\mu_A^c(a),1-\mu_A^c(x)\}\leq\min\{1-\mu_A^c(y),1-\mu_A^c(z)\}.
\]
Hence
\[
\max\{\mu_A^c(y),\mu_A^c(z)\}\leq\max\{\mu_A^c(a),\mu_A^c(x)\},
\]
and the fourth condition of Definition \ref{def21} is satisfied.
\end{proof}

\begin{lemma}\label{lem23} If $A=(\mu_A, \lambda_A)$ is an
$IFSH$ of $G$, then so is \ $\Diamond A=(\lambda^c_A,\lambda_A)$.
\end{lemma}

\begin{proof} The proof is similar to the proof of Lemma
\ref{lem22}. \end{proof}

Combining the above two lemmas it is not difficult to see that the
following theorem is valid.
\begin{theorem}\label{th24} $A=(\mu_A, \lambda_A)$ is an
$IFSH$ of $\,G$ if and only if \ $\Box A$ and $\Diamond A$ are
$IFSHs$ of $\,G$. \hfill $\Box$
\end{theorem}
\begin{corollary}\label{cor25} $A=(\mu_A, \lambda_A)$ is an
$IFSH$ of $\,G$ if and only if \ $\mu_A$ and $\lambda^c_A$ are
fuzzy sub-hyperquasigroups of $\,G$. \hfill $\Box$
\end{corollary}

\begin{theorem}\label{th26} If $\,A=(\mu_A, \lambda_A)$ is an
$IFSH$ of $G$ then the upper $t$-level cut $U(\mu_A;t)$ of $\mu_A$
and the lower $t$-level cut $L(\lambda_A;t)$ of $\lambda_A$ are
sub-hyperquasigroups of $G$ for every \ $t\in Im (\mu_A)\cap
Im(\lambda_A)$.
\end{theorem}

\begin{proof}
Let \ $t\in Im (\mu_A)\cap Im(\lambda_A) \subseteq [0,1]$ and let
$x,y\in  U(\mu_A; t)$. Then $\mu_A(x)\geq t$ and $\mu_A(y)\geq t$
and so $\min\{\mu_A (x),\mu_A(y)\}\geq t$. It follows from the
first condition of Definition \ref{def21} that $\inf\{\mu_A(z)\ |
\ z\in x\circ y\}\geq t$. Therefore for all $z\in x\circ y$ we
have $z\in U(\mu_A;t)$, so $x\circ y\subseteq U(\mu_A;t)$. Hence
for all $a\in U(\mu_A;t)$ we have $a\circ U(\mu_A;t)\subseteq
U(\mu_A;t)$ and $U(\mu_A;t)\circ a\subseteq U(\mu_A;t)$. Now, let
$x\in U(\mu_A;t)$ then there exist $y,z\in G$ such that $x\in
a\circ y$, $x\in z\circ a$ and $\min\{\mu_A(x),\mu_A(a)\}\leq
\min\{\mu(y),\mu(z)\}$. Since $x,a\in U(\mu_A; t)$, we have
$t\leq\min\{\mu_A(x),\mu_A(a)\}$ and so $t\leq\min\{\mu_A(y),
\mu_A(z)\}$ which implies $y\in U(\mu_A;t)$, $z\in U(\mu_A;t)$ and
these prove that $U(\mu_A;t)\subseteq a\circ U(\mu_A;t)$ and
$U(\mu_A;t)\subseteq U(\mu_A;t)\circ a$. Hence $a\circ U(\mu_A;t)=
U(\mu_A;t)= U(\mu_A;t)\circ a$.

Now let $x,y\in L(\lambda_A;t)$. Then $\lambda_A(x)\leq t$, \
$\lambda_A (y)\leq t$ and, consequently, $\max \{\lambda_A
(x),\lambda_A(y)\} \leq t$. It follows from the third condition of
Definition \ref{def21} that $\sup\{\lambda_A(z)\ |\ z\in x\circ
y\}\leq t$. Therefore for all $z\in x\circ y$ we have $z\in
L(\lambda_A;t)$, so $x\circ y \subseteq L(\lambda_A;t)$. Hence for
all $a\in L(\lambda_A;t)$ we have $a\circ L(\lambda_A;t)\subseteq
L(\lambda_A;t)$ and $L(\lambda_A;t)\circ a \subseteq
L(\lambda_A;t)$. Now, let $x\in L(\lambda_A;t)$. Then there exist
$y,z\in G$ such that $x\in a\circ y$, $x\in z\circ a $ and
$\max\{\lambda_A(y),\lambda_A(z)\}\leq
\max\{\lambda(a),\lambda(x)\}$. Since $x,a\in L(\lambda_A;t)$, we
have $\max\{\lambda_A(a), \lambda_A(x)\}\leq t$ and so
$\max\{\lambda_A(y),\lambda_A(z)\}\leq t$ which implies $y\in
L(\lambda_A;t)$, $\,z\in L(\lambda_A;t)$ and these prove that
$L(\lambda_A;t)\subseteq a\circ L(\lambda_A;t)$ and
$L(\lambda_A;t)\subseteq L(\lambda_A;t)\circ a$. Thus $a\circ
L(\lambda_A;t) = L(\lambda_A;t) = L(\lambda_A;t)\circ a$.
\end{proof}

\begin{theorem}\label{th27} If $A=(\mu_A, \lambda_A)$ is an
intuitionistic fuzzy set in $G$ such that the non-empty sets
$U(\mu_A;t)$ and $L(\lambda_A;t)$ are sub-hyperquasigroups of $G$
for all \ $t\in [0,1],$ then $A=(\mu_A, \lambda_A)$ is an $IFSH$
of $G$.
\end{theorem}

\begin{proof} For $t\in [0,1]$, assume that $U(\mu_A;t)\neq
\emptyset$ and $L(\lambda_A;t)\neq\emptyset$ are
sub-hyperquasigroups of $G$. We must show that $A=(\mu_A,
\lambda_A)$ satisfies the all conditions in Definition
\ref{def21}. Let $x,y\in G$, we put $t_0=\min\{\mu_A
(x),\mu_A(y)\}$ and $t_1=\max\{\lambda_A(x),\lambda_A(y)\}$. Then
$x,y\in U(\mu_A;t_0)$ and $x,y\in L(\lambda_A;t_1)$. So $x\circ y
\subseteq U(\mu_A;t_0)$ and $x\circ y \subseteq L(\lambda_A;t_1)$.
Therefore for all $z\in x\circ y$ we have $\mu_A(z)\geq t_0$ and
$\lambda_A(z)\leq t_1$ which imply
\[
\inf\{\mu_A(z) \ | \ z\in x\circ y\}\geq\min\{\mu_A(x),\mu_A(y)\}
\]
 and
\[
\sup\{\lambda_A (z) \ | \ z\in x\circ y\}\leq\max\{\lambda_A(x),
\lambda_A(y)\}
\]
The conditions $(1)$ and $(3)$ of Definition \ref{def21} are
verified.

Now, let $x,a\in G$. If $t_2=\min\{\mu_A(a),\mu_A(x)\}$, then
$a,x\in U(\mu_A;t_2)$. So there exist $y_1,z_1\in U(\mu_A;t_2)$
such that $x\in a\circ y_1$ and $x\in z_1\circ a$. Also we have
$t_2\leq\min\{\mu_A(y_1),\mu_A(z_1)\}$. Therefore the condition
$(2)$ of Definition \ref{def21} is verified. If we put
$t_3=\max\{\lambda_A(a),\lambda_A(x)\}$, then $a,x\in
L(\lambda_A;t_3)$. So there exist $y_2,z_2\in L(\lambda_A; t_3)$
such that $x\in a\circ y_2$ and $x\in z_2\circ a$ and we have
$\max\{\lambda_A(y_2),\lambda_A(y_2)\}\leq t_3$, and so the
condition $(4)$ of Definition \ref{def21} is verified. This
completes the proof.
 \end{proof}

\begin{corollary}\label{cor28} Let $K$ be a sub-hyperquasigroup of
a hyperquasigroup $(G,\circ)$. If fuzzy sets $\mu$ and $\lambda$
are defined on $G$ by
\[
\mu(x)=\left\{\begin{array}{cl}\alpha_0 &\text{if \, $x\in K$},\\
               \alpha_1 &\text{if \, $x\in G\setminus K$,}\end{array}
\right. \qquad\lambda(x)=\left\{\begin{array}{cl}\beta_0
&\text{if \, $x\in K$},\\
                  \beta_1 &\text{if \, $x\in G\setminus K$,}\end{array}
\right.
\]
where $\,0\leq\alpha_1< \alpha_0$, \ $0\leq\beta_0<\beta_1\,$ and
$\,\alpha_i+\beta_i\leq 1$ for $i=0,1,$ then $A=(\mu,\lambda )$ is
an $IFSH$ of $\,G$ and $\,U(\mu;\alpha_0 )=K=L(\lambda;\beta_0)$.
\hfill $\Box$
\end{corollary}

\begin{corollary}\label{cor29}
Let $\chi_{_K}$ be the characteristic function of a
sub-hyperquasigroup $K$ of $\,(G,\circ )$. Then
$K=(\chi_{_K},\chi^c_{_K})$ is an $IFSH$ of $\,G$.
 \hfill $\Box$
\end{corollary}

\begin{theorem}\label{th210} If $\,A=(\mu_{A},\lambda_A)\,$ is an
$IFSH$ of $\,G$, then for all $\,x\in G$ we have
\\[2mm]
\centerline{$\mu_{A}(x)=\sup\{\alpha \in [0,1]\ |\ x\in
U(\mu_{A};\alpha ) \}$}
\\
and\\
\centerline{$ \lambda_A(x)=\inf\{\alpha \in [0,1]\ |\ x\in
L(\lambda_A ;\alpha)\}.$ }
\end{theorem}

\begin{proof} Let $\,\delta=\sup\{\alpha\in [0,1]\ |\ x\in
U(\mu_{A};\alpha )\}\,$ and let $\,\varepsilon >0\,$ be given.
Then $\delta -\varepsilon <\alpha$ for some $\,\alpha \in [0,1]$
such that $\,x\in U(\mu_{A};\alpha)$. This means that $\,\delta
-\varepsilon <\mu_{A}(x)\,$ so that $\,\delta \leq\mu_{A}(x)\,$
since $\,\varepsilon\,$ is arbitrary.

We now show that $\,\mu_{A}(x)\leq\delta.$ If
$\,\mu_{A}(x)=\beta$, then $\,x\in U(\mu_{A};\beta)\,$ and so
\[
\beta\in\{\alpha\in [0,1]\ |\ x\in U(\mu_{A};\alpha )\}.
\]
Hence
\[
\mu_{A}(x)=\beta\leq\sup\{\alpha\in [0,1]\ |\ x\in
U(\mu_{A};\alpha)\}=\delta .
\]
Therefore
\[
\mu_{A}(x)=\delta=\sup\{\alpha\in [0,1]\ |\ x\in
U(\mu_{A};\alpha)\}.
\]
Now let $\,\eta =\inf\{\alpha \in [0,1]\ |\ x\in
L(\lambda_A;\alpha)\}$. Then
\[
\inf\{\alpha\in [0,1]\ |\ x\in L(\lambda_A;\alpha)\}<\eta
+\varepsilon
\]
for any  $\,\varepsilon >0,\,$ and so $\,\alpha <\eta
+\varepsilon\,$ for some $\,\alpha\in [0,1]\,$ with $\,x\in
L(\lambda_A ;\alpha)$. Since $\lambda_A(x)\leq\alpha\,$ and
$\,\varepsilon\,$ is arbitrary, it follows that
$\,\lambda_A(x)\leq\eta$.

To prove $\,\lambda_A(x)\geq\eta$, let $\,\lambda_A(x)=\zeta$.
Then $\,x\in L(\lambda_A ;\zeta)\,$ and thus
$\,\zeta\in\{\alpha\in [0,1]\ |\ x\in L(\lambda_A ;\alpha)\}$.
Hence
\[
\inf\{\alpha\in [0,1]\ |\ x\in L(\lambda_A ;\alpha)\}\leq\zeta,
\]
i.e. $\;\eta\leq\zeta =\lambda_A(x).$ Consequently
\[
\lambda_A(x)=\eta=\inf\{\alpha\in [0,1]\ |\ x\in L(\lambda_A
;\alpha)\},
\]
which completes the proof.
\end{proof}

\begin{theorem}\label{th211} Let $\Omega$ be a non-empty finite subset of
$\,[0,1]$. If $\,\{K_{\alpha}\ |\ \alpha\in\Omega \}$ is a
collection of sub-hyperquasigroups of $\,G\,$ such that
 \begin{itemize}
 \item[\rm (i)] $G= \bigcup\limits_{\alpha\in\Omega} K_{\alpha}$,
 \item[\rm (ii)] $\alpha >\beta\,\Longleftrightarrow\,
        K_{\alpha}\subset K_{\beta}$ \ for all $\alpha,\beta\in\Omega$,
\end{itemize}
then an intuitionistic fuzzy set $A=(\mu_A,\lambda_A)$ defined on
$G$ by

\vspace{4pt}$\mu_{A}(x)=\sup\{\alpha\in\Omega\ |\ x\in
K_{\alpha}\}$ \  \  and \ \ $\lambda_A(x)=\inf\{\alpha\in\Omega\
|\ x\in K_{\alpha}\}$

\vspace{2pt}\noindent is an $IFSH$ of \ $G$.
\end{theorem}
\begin{proof} According to Theorem \ref{th27}, it is
sufficient to show that the non-empty sets $U(\mu_{A};\alpha)$ and
$L(\lambda_A;\beta)$ are sub-hyperquasigroups of $G$. We show that
$U(\mu_{A};\alpha)=K_{\alpha}$. This holds, since
\[
\begin{array}{ll}
x \in U(\mu_{A};\alpha) & \Longleftrightarrow \mu_A(x)\geq \alpha \\
& \Longleftrightarrow \sup \{\gamma \in \Omega \ | \
             x \in K_{\gamma} \} \geq \alpha \\
& \Longleftrightarrow \exists \gamma_0 \in \Omega , \
             x \in K_{\gamma_0}, \ \gamma_0 \geq \alpha \\
& \Longleftrightarrow  x \in K_{\alpha} \ \ \   ({\rm since} \
             K_{\gamma_0} \subseteq K_{\alpha} ).
\end{array}
\]
Now, we prove that $L(\lambda ; \beta )\not = \emptyset$ is a
sub-hyperquasigroup of $G$. We have
\[
\begin{array}{ll}
x \in L(\lambda_{A};\beta) & \Longleftrightarrow \lambda_A(x)\leq \beta \\
& \Longleftrightarrow \inf \{\gamma \in \Omega \ | \
                    x \in K_{\gamma} \} \leq \beta \\
& \Longleftrightarrow \exists \gamma_0 \in \Omega , \
                    x \in K_{\gamma_0}, \ \gamma_0 \leq \beta \\
& \Longleftrightarrow x \in \displaystyle \bigcup_{\gamma \leq \beta} K_{\gamma}
\end{array}
\]
and hence
$L(\lambda_A ; \beta )= \displaystyle \bigcup_{\gamma \leq \beta} K_{\gamma}$.
It is not difficult to see that the union of any family of increasing sub-hyperquasigroups
of a given hyperquasigroup is a sub-hyperquasigroup.
This completes the proof.
\end{proof}

\section{Relations}

Let $\alpha\in [0,1]$ be fixed and let $IFSH(G)$ be the family of
all intuitionistic fuzzy sub-hyperquasigroups of a hyperquasigroup
$G$. For any $A=(\mu_A,\lambda_A)$ and $B=(\mu_B,\lambda_B)$ from
$IFSH(G)$ we define two binary relations $\frak U^{\alpha}$ and
$\frak L^{\alpha}$ on $IFSH(G)$ as follows:
\[
(A,B)\in\frak U^{\alpha}\Longleftrightarrow U(\mu_{A};\alpha) =
U(\mu_B;\alpha)
\]
and
\[
(A,B)\in\frak L^{\alpha}\Longleftrightarrow L(\lambda_A
;\alpha)=L(\lambda_B;\alpha)\, .
 \]
These two relations $\frak U^{\alpha}$ and $\frak L^{\alpha}$ are
equivalence relations. Hence $IFSH(G)$ can be divided into the
equivalence classes of $\frak U^{\alpha}$ and $\frak L^{\alpha}$,
denoted by $[A]_{\frak U^{\alpha}}$ and $[A]_{\frak L^{\alpha}}$
for any $A=(\mu_A,\lambda_A)\in IFSH(G)$, respectively. The
corresponding quotient sets will be denoted by $IFSH(G)/\frak
U^{\alpha}$ and $IFSH(G)/\frak L^{\alpha}$, respectively.

For the family $S(G)$ of all sub-hyperquasigroups of $G$ we define
two maps $U_{\alpha}$ and $L_{\alpha}$ from $IFSH(G)$ to
$S(G)\cup\{\emptyset\}$ by putting
\[
U_{\alpha}(A)=U(\mu_{A};\alpha) \ \ \text{ and } \ \
L_{\alpha}(A)=L(\lambda_A;\alpha)
\]
for each $A=(\mu_A,\lambda_A)\in IFSH(G)$.

It is not difficult to see that these maps are well-defined.

\begin{lemma}\label{lem31} For any $\alpha\in
(0,1)$ the maps $U_{\alpha}$ and $L_{\alpha}$ are surjective.
\end{lemma}
\begin{proof} Let $\bold 0$ and $\bold 1$ be
fuzzy sets in $G$ defined by $\bold 0(x)=0$ and $\bold 1(x)=1$ for
all $x\in G$. Then $\bold 0_{\sim}=(\bold 0,\bold 1)\in IFSH(G)$
and $U_{\alpha}(\bold 0_{\sim})=L_{\alpha}(\bold 0_{\sim})=
\emptyset$ for any $\alpha\in (0,1)$. Moreover for any $K\in S(G)$
we have $K_{\sim} = (\chi_{_K},\chi^c_{_K})\in IFSH(G)$, \
$U_{\alpha}(K_{\sim}) =U(\chi_{_K};\alpha)=K$ and
$L_{\alpha}(K_{\sim})=L(\chi^c_{_K};\alpha)=K$. Hence $U_{\alpha}$
and $L_{\alpha}$ are surjective.
 \end{proof}

\begin{theorem}\label{th31} For any $\alpha\in
(0,1)$ the sets $IFSH(G)/\frak U^{\alpha}$ and $IFSH(G)/\frak
L^{\alpha}$ are equipotent to $S(G)\cup\{\emptyset\}$.
\end{theorem}
\begin{proof} Let $\alpha\in (0,1)$. Putting
$U^*_{\alpha}([A]_{\frak U^{\alpha}}) = U_{\alpha}(A)$ and
$L^*_{\alpha}([A]_{\frak L^{\alpha}})=L_{\alpha}(A)$ for any
$A=(\mu_A,\lambda_A)\in IFSH(G)$, we obtain two maps
\[
U^*_{\alpha}:IFSH(G)/\frak U^{\alpha}\to  S(G)\cup \{\emptyset\} \
\ {\rm and } \ \ L^*_{\alpha}:IFSH(G)/\frak L^{\alpha}\to
S(G)\cup\{\emptyset\}.
 \]
If $U(\mu_{A};\alpha)=U(\mu_B;\alpha)$ and
$L(\lambda_A;\alpha)=L(\lambda_B;\alpha)$ for some
$A=(\mu_A,\lambda_A)$ and $B=(\mu_B,\lambda_B)$ from $IFSH(G)$,
then $(A,B)\in\frak U^{\alpha}$ and $(A,B)\in\frak L^{\alpha}$,
whence $[A]_{\frak U^{\alpha}}=[B]_{\frak U^{\alpha}}$ and
$[A]_{\frak L^{\alpha}}=[B]_{\frak L^{\alpha}}$, which means that
$U*_{\alpha}$ and $L^*_{\alpha}$ are injective.

To show that the maps $U^*_{\alpha}$ and $L_{\alpha}$ are
surjective, let $K\in S(G)$. Then for
$K_{\sim}=(\chi_{_K},\chi^c_{_K})\in IFSH(G)$ we have
$U^*_{\alpha}([K_{\sim}]_{\frak U^{\alpha}}) =
U(\chi_{_K};\alpha)=K$ and $L^*_{\alpha}([K_{\sim}]_{\frak
L^{\alpha}})= L(\chi^c_{_K};\alpha)= K$. Also $\bold
0_{\sim}=(\bold 0,\bold 1)\in IFSH(G)$. Moreover
$U^*_{\alpha}([\bold 0_{\sim}]_{\frak U^{\alpha}})=U(\bold
0;\alpha)=\emptyset$ and $L^*_{\alpha}([\bold 0_{\sim}]_{\frak
L^{\alpha}}) = L(\bold 1;\alpha)=\emptyset .$ Hence $U^*_{\alpha}$
and $L^*_{\alpha}$ are surjective.
 \end{proof}

Now for any $\alpha \in [0,1]$ we define a new relation $\frak
R^{\alpha}$ on $IFSH(G)$ by putting:
\[
(A,B)\in\frak R^{\alpha}\Longleftrightarrow U(\mu_{A};\alpha)\cap
L(\lambda_A;\alpha) =U(\mu_B;\alpha)\cap L(\lambda_B;\alpha),
 \]
where $A=(\mu_A,\lambda_A)$ and $B=(\mu_B,\lambda_B)$. Obviously
$\frak R^{\alpha}$ is an equivalence relation.

\begin{lemma}\label{lem33}
The map $\;I_{\alpha}: IFSH(G)\to S(G)\cup\{\emptyset\}$ defined
by
 \[
I_{\alpha}(A)=U(\mu_{A};\alpha)\cap L(\lambda_A;\alpha),
 \]
where $A=(\mu_{A},\lambda_A)$, is surjective for any $\alpha\in
(0,1)$.
\end{lemma}

\begin{proof} If $\alpha\in (0,1)$ is fixed, then for
$\bold 0_{\sim}=(\bold 0,\bold 1)\in IFSH(G)$ we have
 \[
I_{\alpha}(\bold 0_{\sim})=U(\bold 0;\alpha)\cap L(\bold 1;\alpha)
= \emptyset\, ,
 \]
and for any $K\in S(G)$ there exists
$K_{\sim}=(\chi_{_K},\chi^c_{_K})\in IFSH(G)$ such that
$I_{\alpha}(K_{\sim})=U(\chi_{_K};\alpha)\cap
L(\chi^c_{_K};\alpha) = K$.
\end{proof}

\begin{theorem}\label{th34} For any $\alpha\in (0,1)$ the
quotient set $IFSH(G)/\frak R^{\alpha}$ is equipotent to
$S(G)\cup\{\emptyset\}$.
\end{theorem}
\begin{proof} Let $I^*_{\alpha} : IFSH(G)/\frak R^{\alpha}\to
S(G)\cup\{\emptyset\}$, where $\alpha\in (0,1)$, be defined by the
formula:
 \[
I^*_{\alpha}([A]_{\frak R^{\alpha}}) = I_{\alpha}(A) \ \ \text{
for each } \ \ [A]_{\frak R^{\alpha}}\in IFSH(G)/\frak R^{\alpha}.
 \]
If $I^*_{\alpha}([A]_{\frak R^{\alpha}})=I^*_{\alpha}([B]_{\frak
R^{\alpha}})$ for some $[A]_{\frak R^{\alpha}},\, [B]_{\frak
R^{\alpha}}\in IFSH(G)/\frak R^{\alpha}$,  then
 \[
U(\mu_{A};\alpha)\cap L(\lambda_A;\alpha)=U(\mu_B;\alpha)\cap
L(\lambda_B;\alpha) ,
 \]
which implies $(A,B)\in\frak R^{\alpha}$ and, in the consequence,
$[A]_{\frak R^{\alpha}}=[B]_{\frak R^{\alpha}} .$ Thus
$I^*_{\alpha}$ is injective.

It is also onto because $I^*_{\alpha}(\bold
0_{\sim})=I_{\alpha}(\bold 0_{\sim})=\emptyset$ for $\bold
0_{\sim}=(\bold 0,\bold 1)\in IFSH(G)$, and
$I^*_{\alpha}(K_{\sim}) = I_{\alpha}(K) = K$ for $K\in S(G)$ and
$K_{\sim}=(\chi_{_K},\chi^c_{_K})\in IFSH(G)$.
\end{proof}

\section{Connections with binary quasigroups}

A groupoid $(Q,\cdot)$ is called a (\emph{binary})
\emph{quasigroup} if each of the equations $ax=b$ and $ya=b$ has a
unique solution for any $a,b\in Q$. Since a non-empty subset of
$Q$ closed with respect to this operation is not in general a
quasigroup we must use the another equivalent definition of a
quasigroup. A quasigroup $(Q,\cdot)$ can be defined (cf.
\cite{olga}) as an algebra $(Q,\cdot,\setminus, \, / )$ with three
binary operation such that $(Q,\cdot )$ is a quasigroup in the
above sense and
\[
x\setminus y=z\Leftrightarrow xz=y \ \ \text{ and } \ \ x/ y=z
\Leftrightarrow zy=x
\]
for all $x,y,z\in Q$. In this case a non-empty subset of $Q$ is a
subquasigroup of $(Q,\cdot)$ (and $(Q,\cdot ,\setminus,\, /)$) if
and only if it is closed with respect to these three operations.
This gives the possibility to the introduction of a good
definition of intuitionistic fuzzy subquasigroups of binary
quasigroups \cite{11}.

\begin{definition}\label{def35} Let $(Q,\cdot)$ be a quasigroup.
An intuitionistic fuzzy set $A=(\mu_A, \lambda_A)$ in $Q$ is
called an \emph{intuitionistic fuzzy subquasigroup} of $Q$ if
\begin{itemize}
\item[(i)] $\min \{ \mu_A(x), \mu_A(y) \} \leq \mu_A(x*y)$
\item[(ii)] $\lambda_A(x*y)\leq \max \{ \lambda_A(x), \lambda_A(y)
\}$
\end{itemize}
hold for all $x,y\in Q$ and $*\in\{\cdot,\setminus,\, /\}.$
\end{definition}

In this case an intuitionistic fuzzy set $A=(\mu_A, \lambda_A)$ is
an intuitionistic fuzzy subquasigroup of $(Q,\cdot,\setminus, \, /
)$ if and only if all non-empty $U(\mu;t)$ and $L(\mu;t)$ are
subquasigroups of $(Q,\cdot,\setminus, \, / )$ (cf. \cite{11}).

A hyperquasigroup $(G, \circ )$ is called \emph{regular} if
\[
x \in y\circ z \ \ {\rm implies } \ \ y\in x \circ z \ \ {\rm and
} \ \ z \in y \circ x
\]
for all $x,y,z \in G$. Let $(G,\circ )$ be a regular
hyperquasigroup. The relation $\beta^*$ is the smallest
equivalence relation on $G$ such that the quotient $G/\beta^*$,
the set of all equivalence classes, is a quasigroup. $\beta^*$ is
called the fundamental equivalence relation on $G$ and $G/\beta^*$
is called the fundamental quasigroup. The equivalence relation
$\beta^*$ was introduced by Koskas \cite{13} and studied mainly by
Corsini \cite{4} and Freni \cite{9}, \cite{10} concerning
hypergroups and Vougiouklis \cite{16} concerning $H_v$-groups.

Let us denote by ${\cal U}$ the set of all finite products of
elements of $G$ as follows:
\[
x \beta y \ {\rm if \ and \ only \ if } \ \{x,y \} \subseteq u \
{\rm for \ some} \ u \in {\cal U}.
\]
The fundamental relation $\beta^*$ is the transitive closure of
the relation $\beta$ (see Theorem 1.2.2 in \cite{16}). Suppose
$\beta^*(a)$ is the equivalence class containing $a \in G$. Then
the product ``$\cdot$'' on $G/\beta^*$ is defined as follows:
\[
\beta^*(a)\cdot\beta^*(b)=\beta^*(c)\ \ {\rm for \ all} \ \
c\in\beta^*(a)\circ\beta^*(b).
\]
In this case, each of the equations $\beta^*(a)\cdot
\beta^*(x)=\beta^*(b)$ and $\beta^*(y)\cdot \beta^*(a)=\beta^*(b)$
has a unique solution for any $\beta^*(a), \beta^*(b) \in
G/\beta^*$. The quasigroup $(G/\beta^*, \cdot , \setminus , / )$
corresponds to quasigroup $(G/\beta^* , \cdot )$, where
\[
\begin{array}{lll}
\beta^*(x)\setminus\beta^*(y)=\beta^*(z) & \Longleftrightarrow &
\beta^*(x) \cdot \beta^*(z)=\beta^*(y),\\
\beta^*(x) / \beta^*(y)=\beta^*(z) & \Longleftrightarrow &
\beta^*(z) \cdot \beta^*(y)=\beta^*(x).
\end{array}
\]
Let $\mu$ be a fuzzy set in $G$. The fuzzy set $\mu_{\beta^*}$ in
$G/\beta^*$ is defined as follows:
\[ \mu_{\beta^*}:G/\beta^*\to  [0,1], \quad
   \beta^*(x)\mapsto \sup\{\mu(a)\ | \ a\in\beta^*(x)\}. \]
Now, we have

\begin{theorem}\label{th36} Let $G$ be a regular hyperquasigroup and
$A=(\mu_A, \lambda_A)$ an intuitionistic fuzzy sub-hyperquasigroup
of $G$. Then $A/\beta^*=(\mu_{\beta^*}, \lambda_{\beta^*})$ is an
intuitionistic fuzzy subquasigroup of the fundamental quasigroup
$G/\beta^*$.
\end{theorem}

\medskip\noindent
{\bf Acknowledgements.} The authors are highly grateful to the
referees for their valuable comments and suggestions for improving
the paper.

\end{document}